\documentclass[11pt,dvipdfm]{article}
\usepackage{amsfonts}
\usepackage{graphicx}
\usepackage{amsfonts, amsmath, amssymb, latexsym, epsfig}
\usepackage{mathrsfs}
\newcommand{\blst}{\begin{trivlist}}
\newcommand{\elst}{\end{trivlist}}

\newtheorem{thm}{Theorem}[section]
\newtheorem{prop}[thm]{Proposition}
\newtheorem{cor}[thm]{Corollary}
\newtheorem{lem}[thm]{Lemma}

\newtheorem{conj}[thm]{Conjecture}
\newtheorem{exa}[thm]{Example}
\newtheorem{defn}[thm]{Definition}

\newcommand{\ben}{\begin{enumerate}}
\newcommand{\een}{\end{enumerate}}
\newcommand{\ble}{\begin{lem}}
\newcommand{\ele}{\end{lem}}
\newcommand{\bth}{\begin{thm}}
\renewcommand{\eth}{\end{thm}}
\newcommand{\bpr}{\begin{prop}}
\newcommand{\epr}{\end{prop}}
\newcommand{\bco}{\begin{cor}}
\newcommand{\eco}{\end{cor}}
\newcommand{\bcon}{\begin{conj}}
\newcommand{\econ}{\end{conj}}
\newcommand{\bde}{\begin{defn}}
\newcommand{\ede}{\end{defn}}
\newcommand{\bex}{\begin{exa}}
\newcommand{\eex}{\end{exa}}
\newcommand{\barr}{\begin{array}}
\newcommand{\earr}{\end{array}}
\newcommand{\btab}{\begin{tabular}}
\newcommand{\etab}{\end{tabular}}
\newcommand{\beq}{\begin{equation}}
\newcommand{\eeq}{\end{equation}}

\newcommand{\bea}{\begin{eqnarray*}}
\newcommand{\eea}{\end{eqnarray*}}
\newcommand{\beaa}{\begin{eqnarray}}
\newcommand{\eeaa}{\end{eqnarray}}
\newcommand{\bce}{\begin{center}}
\newcommand{\ece}{\end{center}}
\newcommand{\bpi}{\begin{picture}}
\newcommand{\epi}{\end{picture}}
\newcommand{\bfi}{\begin{figure} \begin{center}}
\newcommand{\efi}{\end{center} \end{figure}}

\newcommand{\bsl}{\begin{slide}{}}
\newcommand{\esl}{\end{slide}}

{}
{}
{}
{}
{}
{}
{}
{}
{}
{}
{}
{}
{}
{}
{}
{}
{}
{}
{}
{}
{}
{}
{}
{}
{}
{}
{}
{}
\newenvironment{proof}{
\par
\noindent {\bf Proof.}\rm}{\mbox{}\hfill\rule{0.5em}{0.809em}\par}
\begin{document}

\title{Tutte polynomial and $G$-parking functions}
\author{Hungyung Chang$^{a}$
\and
 \and Jun Ma$^{b,}$\thanks{Email address of the corresponding author: majun@math.sinica.edu.tw}\\
 \and Yeong-Nan Yeh$^{c,}$\thanks{Partially supported by NSC 96-2115-M-001-005}
}\date{} \maketitle \vspace*{-1.2cm}\begin{center} \footnotesize
$^{a,b,c}$ Institute of Mathematics, Academia Sinica, Taipei, Taiwan\\

\end{center}

 \vspace*{-0.3cm}
\thispagestyle{empty}
\begin{abstract}
Let $G$ be a connected graph with vertex set $\{0,1,2,\ldots,n\}$.
We allow $G$ to have multiple edges and loops. In this paper, we
give a characterization of external activity by some parameters of
$G$-parking functions. In particular, we give the definition of the
bridge vertex of a $G$-parking function and obtain an expression of
the Tutte polynomial $T_G(x,y)$ of $G$ in terms of $G$-parking
functions.  We find the Tutte polynomial enumerates the $G$-parking
function by the number of the bridge vertices.
\end{abstract}
%Keyword

%main text
%Section 1
%\newpage
\section{Introduction}
J. Riordan \cite{R} define the parking function as follows: $m$
parking spaces are arranged in a line, numbered $1$ to $n$ left to
right; $n$ cars, arriving successively, have initial parking
preferences, $a_i$ for $i$, chosen independently and at random;
$(a_1,\cdots,a_n)$ is called preference function; if space $a_i$ is
occupied, car $i$ moves to the first unoccupied space to the right;
if all the cars can be parked, then  the preference function is
called parking function.

Konheim and Weiss \cite{konhein1966} introduced the conception of
the  parking functions of length $n$ in the study of the linear
probes of random hashing function.  J. Riordan \cite{R} studied  the
parking functions  and derived that the number of parking functions
of length $n$ is $(n+1)^{n-1}$, which coincides with the number of
labeled trees on $n+1$ vertices by Cayley's formula. Several
bijections between the two sets are known (e.g., see
\cite{FR,R,SMP}). Parking functions have been found in connection to
many other combinatorial structures such as acyclic mappings,
polytopes, non-crossing partitions, non-nesting partitions,
hyperplane arrangements, etc. Refer to \cite{F,FR,GK,PS,SRP,SRP2}
for more information.

Parking function $(a_1,\cdots,a_n)$ can be redefined that its
increasing rearrangement $(b_1,\cdots,b_n)$ satisfies $b_i\leq i$.
 Pitman and  Stanley generalized the notion of parking functions
in \cite{PS}. Let ${\bf x}=(x_1,\cdots,x_n)$ be a sequence of
positive integers. The sequence $\alpha=(a_1,\cdots,a_n)$  is called
an ${\bf x}$-parking function if the non-decreasing rearrangement
$(b_1,\cdots,b_n)$ of $\alpha$ satisfies $b_i\leq x_1+\cdots +x_i$
for any $1\leq i\leq n$. Thus, the ordinary parking function is the
case ${\bf x}=(1,\cdots,1)$. By the determinant formula of
Gon\v{c}arove polynomials, Kung and Yan \cite{KY} obtained the
number of ${\bf x}$-parking functions for an arbitrary ${\bf x}$.
See also \cite{Y1,Y2,Y3} for the explicit formulas and properties
for some specified cases of ${\bf x}$.

Recently, Postnikov and Shapiro \cite{postnikov2004} gave a new
generalization, building on work of Cori, Rossin and Salvy
\cite{cori2002}, the $G$-parking functions of a graph. For the
complete graph $G=K_{n+1}$, the defined functions in
\cite{postnikov2004} are exactly the classical parking functions.
Chebikin and Pylyavskyy \cite{Denis2005} established a family of
bijections from the set of $G$-parking functions to the spanning
trees of $G$.

Dimitrije Kostic and Catherine H. Yan \cite{kostic} proposed the
notion of a $G$-multiparking function, a natural extension of the
notion of a $G$-parking function and extended the result of
\cite{Y3} to arbitrary graphs. They constructed a family of
bijections from the set of $G$-multiparking functions to the
spanning forests of $G$. Particularly, They characterize the
external activity by the bijection induced by the breadth-first
search and gave a representation of Tutte polynomial by the reversed
sum of $G$-multiparking functions. Given a classical parking
function $\alpha=(a_1,\ldots,a_n)$, let $cr(\alpha)$ be the number
of critical left-to-right maxima in $\alpha$. They also gave an
expression of the Tutte polynomial $T_{K_{n+1}}(x,y)$ of the
complete graph $K_{n+1}$ as follows:
$$T_{K_{n+1}}(x,y)=\sum\limits_{\alpha\in\mathcal{P}_n}x^{cr(\alpha)}y^{{n\choose
2}-\sum\limits_{i=1}^na_i},$$ where $\mathcal{P}_n$ is the set of
classical parking functions of length $n$. Recently, Sen-peng Eu,
Tung-Shan Fu and Chun-Ju Lai \cite{EFL} considered a class of
multigraphs in connection with ${\bf x}$-parking functions, where
${\bf x}=(a,b,\ldots,b)$. They gave the Tutte polynomial of the
multigraphs in terms of ${\bf x}$-parking functions.

Let $G$ be a connected graph with vertex set $\{0,1,2,\ldots,n\}$.
We allow $G$ to have multiple edges and loops. The motivation of
this paper is to extend the results in \cite{Y3} on $K_{n+1}$ to
arbitrary connected graphs and give a characterization of external
activity by some parameters of $G$-parking functions. To obtain the
characterization for the complete graph $K_{n+1}$,  Dimitrije Kostic
and Catherine H. Yan \cite{kostic} use the bijections induced by the
breadth-first search. In this paper, we use the bijections induced
by the vertex ranking. We give the definition of the bridge vertex
of a $G$-parking functions. We obtain a expression of the Tutte
polynomial $T_G(x,y)$ of $G$ in terms of $G$-parking functions. So,
we find the Tutte polynomial enumerates the $G$-parking function by
the number of the bridge vertices.

This paper is organized as follows. In Section $2$, we give the
definition of the bridge vertex of a $G$-parking functions. In
Section $3$,  we will express the Tutte polynomial $T_G(x,y)$ of $G$
in terms of $G$-parking functions.
\section{The bridge vertex of $G$-parking function}
In this section, we always let $G$ be a connected graph with vertex
set $\{0,1,2,\ldots,n\}$ and edge set $E(G)$. We allow $G$ to have
multiple edges and loops. Let $[n]:=\{1,2,\ldots,n\}$. For any
$I\subseteq V(G)\setminus \{0\}$ and $v\in I$, define
$outdeg_{I,G}(v)$ to be the cardinality of the set $\{\{w,v\}\in
E(G)\mid w\notin I\}$. We give the definition of $G$-parking
function as follows.
\begin{defn}\label{definition} Let $G$ be a connected graph with vertex set
$V(G)=\{0,1,2,\cdots,n\}$ and edge set $E(G)$. A $G$-parking
function is a function $f:V(G)\rightarrow \mathbb{N}\cup\{-1\}$,
such that for every $I\subseteq V(G)\setminus\{0\}$ there exists a
vertex $v\in I$ such that $0\leq f(v)<outdeg_{I,G}(v)$ and
$f(0)=-1$.
\end{defn}

For any $i,j\in[n]$, let $\mu_G(i,j)$ be the number of edges
connecting the vertices $i$ to $j$ in $G$. For establishing the
bijections, all edges of $G$ are colored. The colors of edges
connecting the vertices $i$ to $j$ are $0,1,\cdots,\mu_G(i,j)-1$
respectively for any $i,j\in V(G)$. We use $\{i,j\}_k$ to denote the
edge $e\in E(G)$ connecting two vertices $i$ and $j$ with color $k$.
A subgraph $T$ of $G$ is called a subtree of $G$ rooted at $m$ if
the subgraph contains the vertex $m$ and there is a unique path from
$i$ to $m$ in $T$ for every vertex $i$ of $T$. If a subtree contains
all vertices of $G$, then we say the subtree is a spanning tree of
$G$. Let $\mathcal{P}_{G}$ and $\mathcal{T}_{G}$ be the sets of the
$G$-parking functions and the color spanning trees of $G$
respectively. For any $T\in\mathcal{T}_{G}$ and $e\in T$, let
$c_T(e)$ denote the color of edge $e$ in $T$. Kostic and Yan
\cite{kostic} give an algorithm $\Phi$ which
is a bijection from the sets $\mathcal{P}_G$ to $\mathcal{T}_G$. We give a description
 of the algorithm as follows.\\

\noindent{\bf Algorithm A. (Kostic, Yan \cite{kostic})}

{\bf Step 1:} Let $val_0=f$, $P_0=\emptyset$, $T_0=Q_0=\{0\}$.

{\bf Step 2:} At time $i\geq 1$, let $v=\min \{\tau(w)\mid w\in
Q_{i-1}\}$, where $\tau$ is a vertex ranking in $S_n$.

{\bf Step 3:} Let $N=\{w\notin P_{i-1}\mid 0\leq val_{i-1}(w)\leq
\mu(w,v)-1\text{ and }\{w,v\}_{val_{i-1}(w)}\in E(G)\}$ and
$\hat{N}=\{w\notin P_{i-1}\mid val_{i-1}(w)\geq \mu(w,v)\text{ and
}\{w,v\}_{val_{i-1}(w)}\in E(G)\}$. Set
$val_i(w)=val_{i-1}(w)-\mu_G(w,v)$ for all $w\in \hat{N}$. For any
other vertex $u$, set $val_i(w)=val_{i-1}(w)$. Update $P_i$, $Q_i$
and $F_i$ by letting $P_i=P_{i-1}\cup\{v\}$, $Q_i=Q_{i-1}\cup
N\setminus \{v\}$. Let $T_i$ be a graph on $P_{i}\cup Q_i$ whose
edges are obtained from those of $T_{i-1}$ by joining edges
$\{w,v\}_{val_{i-1}(w)}$ for each $w\in N$.\\

Define $\Phi=\Phi_{G,\tau}:\mathcal{P}_{G}\rightarrow
\mathcal{T}_{G,}$ by letting $\Phi(f)=T_n$. Let $T\in
\mathcal{F}_{G}$. Note that the vertex $0$ is the root of $T$. For
any non-root vertex $v\in [n]$, there is a unique path from $v$ to
$0$ in $T$. Define the height of $v$ to be the number of edges in
the path. If the height of a vertex $w$ is less than the height of
$v$ and $\{v,w\}_k$ is an edge of $T$, then $w$ is the predecessor
of $v$, $v$ is a child of $w$, and write $w={ pre}_T(v)$ and $v\in {
child}_T(w)$.  Suppose $T'$ is a subtree of $G$. A leaf of $T'$ is a
vertex of $T'$ with degree $1$ in $T'$. Denote the set of leaves of
$T'$ by $Leaf(T')$. The following algorithm will give the inverse
map of
$\Phi$.\\

\noindent{\bf Algorithm B (Kostic, Yan \cite{kostic}).}

{\bf Step 1. } Let $\tau$ be a vertex ranking in $S_n$. Assume
$v_0,v_1,v_2,\ldots, v_i$ are determined, where $v_0=0$. Let
$V_i=\{v_0,v_1,v_2,\ldots, v_i\}$ and $W_i=\{v\notin V_i\mid
\{v,w\}_k\in T\text{ for some }w\in V_i\}$. Let $T'$ be the subtree
obtained by restricting $T$ to $V_i\cup W_i$. Let $v_{i+1}$ be the
vertex $w$ such that $\tau(w)\leq\tau(u)$ for all $u\in Leaf(T')$.

{\bf Step 2. } Let $\pi=(v_1\ldots v_n)$ be the order of the
vertices of $G$ determined by Step 1. Set $f(0)=-1$. For any other
vertex $v$, let $f(v)$ be equal to the sum of
 the color of the edge connecting the vertices $v$ and $pre_T(v)$ and the cardinality of the set
 $N(v)$, where
$N(v)=N_{G,T,\tau}(v)=\{w\mid \{v,w\}_k\in E(G)\text{ and }
\pi^{-1}(w)<\pi^{-1}({\rm pre}_T(v))\}$.\\

Define $\Theta=\Theta_{G,\tau}:\mathcal{T}_{G}\rightarrow
\mathcal{P}_{G}$ by letting $\Theta_{G,\tau}(T)=f_T$. Then $\Theta$
is the inverse of $\Phi$. Note that the order $\pi=v_1v_2\ldots v_n$
in the algorithm B is exactly the order in which vertices of $G$
will be placed into the set $P_i$ when running algorithm A on $f$.
Define $Ord=Ord_{G,\tau}:\mathcal{P}_{G}\rightarrow S_{n}$ by
letting $Ord_{G,\tau}(f)=(v_1v_2\ldots v_n)$, where the order
$0=v_0,v_1,v_2,\ldots v_n$ is obtained by Algorithm B,i.e.,
$Ord_{G,\tau}(f)_i=u$ and $Ord^{-1}_{G,\tau}(f)_u=i$ if $v_i=u$ for
all $i\in [n]$. Furthermore, let $Rea(f)$ be a function such that
$Rea(f)(i)=f(Ord(f)_i)$ for all $i\in [n]$. We say $Ord(f)$ and
$Rea(f)$ are an order and  a rearrangement of the $G$-parking
function $f$ respectively. Hence, for any $f\in\mathcal{P}_G$, we
can obtain a pair $(Rea(f),Ord(f))$.

\begin{defn} Let $f\in\mathcal{P}_G$ and $v\in V(G)$. Suppose $Ord(f)_i=v$. let
$I_v=I_{G,\tau,f,v}=\{Ord(f)_j\mid j\geq i\}$. The vertex $v$ is
said to be {\it $f$-critical} if $f(v)={\rm outdeg}_{I_v}(v)-1$.
\end{defn}

Define $C_f=C_{G,\tau,f}$ to be the set of all the $f$-critical
vertices. Clearly, $C_f\neq\emptyset$ for any $f\in \mathcal{P}_G$
since $0\in C_f$.

\begin{exa} Let us consider the following graph $G$. Let $\tau$ be the identity permutation.\\
\begin{center}\includegraphics[width=3cm]{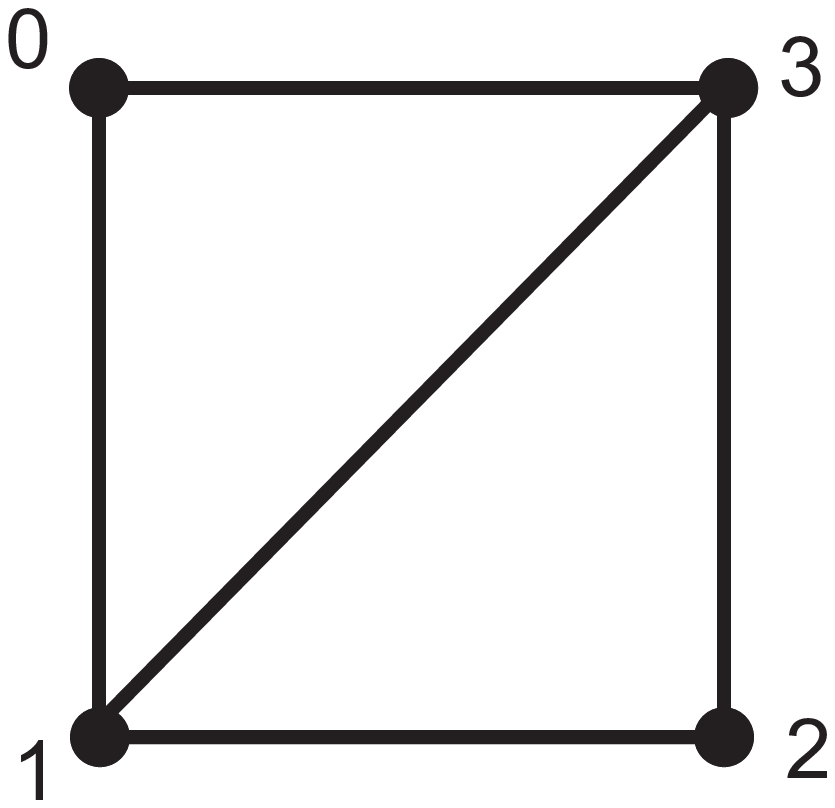}\\
Fig 1. A graph $G$\end{center} We list all the $G$-parking functions
$f$ as well as the corresponding $Ord(f)$, $Rea(f)$ and $C_f$ as
follows.
$$\begin{array}{|l|l|l|l|l|l|}
\hline G-parking~functions~f&Ord(f)&Rea{f}&C_f\\
\hline f_1=(-1,0,0,0)&(0,1,2,3)&(-1,0,0,0)&\{0,1,2\}\\
\hline f_2=(-1,0,0,1)&(0,1,2,3)&(-1,0,0,1)&\{0,1,2\}\\
\hline f_3=(-1,0,0,2)&(0,1,2,3)&(-1,0,0,2)&\{0,1,2,3\}\\
\hline f_4=(-1,0,1,0)&(0,1,3,2)&(-1,0,0,1)&\{0,1,2\}\\
\hline f_5=(-1,0,1,1)&(0,1,3,2)&(-1,0,1,1)&\{0,1,2,3\}\\
\hline f_6=(-1,1,0,0)&(0,3,1,2)&(-1,0,1,0)&\{0,1,3\}\\
\hline f_7=(-1,1,1,0)&(0,3,1,2)&(-1,0,1,1)&\{0,1,2,3\}\\
\hline f_8=(-1,2,0,0)&(0,3,2,1)&(-1,0,0,2)&\{0,1,2,3\}\\
\hline
\end{array}
$$
\begin{center}Table 1. All the $G$-parking functions
\end{center}
\end{exa}

\begin{defn} Let $f\in\mathcal{P}_G$ and $v\in V(G)\setminus \{0\}$.
Suppose $Ord(f)_i=v$. A $G$-parking function $g$ is {\it weak
$v$-identical to $f$} if it satisfies the following conditions:

(1) $Rea(g)(j)=Rea(f)(j)$ and $Ord(g)_j=Ord(f)_j$ for all $j\in
[i-1]$,

(2) $g(v)\geq f(v)$, and

(3) $g(w)\geq {\rm outdeg}_{I_{v}}(w)$
 for all $w\in I_v$ and $\tau(w)<\tau(v)$.\\
Furthermore, $g$ is   {\it strong $v$-identical to $f$} if (1) $g$
is weak $v$-identical to $f$; (2) $Ord(g)_i=v$.
\end{defn}

Given $f\in\mathcal{P}_G$ and $v\in V(G)\setminus \{0\}$, define
$$W_{v,f}=W_{G,v,\tau,f}=\{g\in\mathcal{P}_G\mid g\text{ is weak
}v\text{-identical to }f\}$$ and
$$S_{v,f}=S_{G,v,\tau,f}=\{g\in\mathcal{P}_G \mid g\text{ is strong }\text{
}v\text{-identical to }f\}.$$ It is easy to see that
$S_{v,f}\subseteq W_{v,f}$ and $g(v)=f(v)$ for all $g\in S_{v,f}$ if
$v\in C_f$.

\begin{lem}\label{nobridge} Let $G$ be a connected graph and $f$ a $G$-parking function.
Let $e$ be an edge of $G$ connecting the vertices $w$ to $v$.
Suppose that $e$ is a bridge of $G$ and the vertices $w$ and $0$ are
in the same component after deleting the edge $e$. Then $v\in B(f)$
for any $f\in\mathcal{P}_G$.
\end{lem}
\begin{proof} Since $e$ is a bridge of $G$ and the vertices $w$ and $0$ are
in the same component after deleting the edge $e$, we have $f(v)=0$,
$Ord(f)^{-1}(w)<Ord(f)^{-1}(v)$ and $v\in C_f$ for all
$f\in\mathcal{P}_G$. Given a $f\in\mathcal{P}_G$, suppose
$Ord(f)_i=v$. Assume that $W_{v,f}\neq S_{v,f}$, i.e., there is a
$g\in\mathcal{P}_G$ such that $g\in W_{v,f}$ and $g\notin S_{v,f}$.
$g\notin S_{v,f}$ implies $Ord(g)_i\neq v$. Let $u=Ord(g)_i$. Then
$\tau(u)>\tau(v)$ and $Ord(g)^{-1}(u)<Ord(g)^{-1}(v)$ since $g\in
W_{v,f}$. Note that
$Ord(g)^{-1}(w)=Ord(f)^{-1}(w)<Ord(f)^{-1}(v)=Ord(g)^{-1}(u)$. So,
we must have $Ord(g)^{-1}(v)<Ord(g)^{-1}(u)$ by Algorithm A, a
contradiction. Hence, $W_{v,f}= S_{v,f}$ and $v\in B(f)$.
\end{proof}

\begin{defn} Let $f\in\mathcal{P}_G$. A vertex $v\in V(G)\setminus\{0\}$ is said to be $f$-bridge if
 $v\in C_f$ and $|W_{v,f}|=|S_{v,f}|$.
\end{defn}
 Define $B(f)=B_{G,\tau}(f)$ as the set of the $f$-bridge
 vertices of $f$, $b(f)=b_{G,\tau}(f)=|B_{G,\tau}(f)|$ and
$w(f)=w_{G}(f)=|E(G)|-|V(G)|-\sum\limits_{i=0}^nf(i)$.

\begin{exa} We consider the graph $G$ in Fig 1. By Table 1, we have
$f_3=(-1,0,0,2)$ is a $G$-parking function and
$C_{f_3}=\{0,1,2,3\}$. It is easy to check the results in Table 2.
$$\begin{array}{|l|l|l|l|l|l|}
\hline W_{1,f_3}=\{f_i\mid i\in [8]\}&S_{1,f_3}=\{f_i\mid i\in [5]\}\\
\hline W_{2,f_3}=\{f_i\mid i\in [5]\}&S_{2,f_3}=\{f_i\mid i\in [3]\}\\
\hline W_{3,f_3}=\{f_i\mid i\in [3]\}&S_{3,f_3}=\{f_i\mid i\in
[3]\}\\
\hline
\end{array}
$$
\begin{center} Table 2. $W_{v,f_3}$ and $S_{v,f_3}$, where $v\in C_{f_3}$
\end{center} Hence, $B(f_3)=\{3\}$. We list all the $G$-parking functions as
well as the corresponding parameters $b(f)$ and $w(f)$ in the
following table.
$$\begin{array}{|l|l|l|}
\hline G-parking function~f&B(f)&(b(f),w(f))\\
\hline f_1=(-1,0,0,0)&\emptyset&(0,2)\\
\hline f_2=(-1,0,0,1)&\emptyset&(0,1)\\
\hline f_3=(-1,0,0,2)&\{3\}&(1,0)\\
\hline f_4=(-1,0,1,0)&\{2\}&(1,1)\\
\hline f_5=(-1,0,1,1)&\{2,3\}&(2,0)\\
\hline f_6=(-1,1,0,0)&\{3\}&(1,1)\\
\hline f_7=(-1,1,1,0)&\{2,3\}&(2,0)\\
\hline f_8=(-1,2,0,0)&\{1,2,3\}&(3,0)\\
\hline
\end{array}
$$
\begin{center}
Table 3. $G$-parking functions $f$ as well as the corresponding
parameters $b(f)$ and $w(f)$
\end{center}
We note that the Tutte polynomial $T_G(x,y)$ of $G$ in Fig.1
satisfies
$$T_G(x,y)=x^3+2x^2+x+2xy+y+y^2=\sum\limits_{f\in\mathcal{P}_G}x^{b(f)}y^{w(f)}.$$
\end{exa}

\section{A new expression
of the Tutte polynomial} In this section, we will prove the main
theorem of this paper. Suppose that $e$ is an edge connecting the
vertices $i$ to $j$ in $G$, where $i<j$. Define a graph $G{\setminus
e}$ as follows. The graph $G{\setminus e}$ is obtained from $G$
contracting the the vertices $i$ and $j$; that is, to get
$G{\setminus e}$ we identify two vertices $i$ and $j$ as a new
vertex $i$. Define $G-e$ as a graph obtained by deleting the edge
$e$ from $G$.

Let $NB_G(0)$ be the set of the vertices which are adjacent to the
vertex $0$ in $G$. We consider the case in which $e$ is an edge of
$G$ connecting the root $0$ to the vertex $u$ and $u$ satisfies
$\tau(u)\leq\tau(w)$ for all $w\in NB_G(0)$. Let
$\mathcal{P}_G^0=\{f\in\mathcal{P}_G\mid f(u)=0\}$ and
$\mathcal{P}_G^1=\{f\in\mathcal{P}_G\mid f(u)\geq 1\}$. Clearly,
$\mathcal{P}_G^0\cap\mathcal{P}_G^1=\emptyset$ and
$\mathcal{P}_G=\mathcal{P}_G^0\cup\mathcal{P}_G^1$. For any $f\in
\mathcal{P}_G^0$, let $g=\phi(f)$ such that $g(w)=f(w)$ for any
$w\neq u$. For any $f\in \mathcal{P}_G^1$, let $g=\varphi(f)$ such
that $g(w)=f(w)$ for any $w\neq u$ and $g(u)=f(u)-1$.
\begin{lem}\label{lemmabijectionu=0} (1) The mapping $\phi$ is a bijection from $\mathcal{P}^0_G$
to $\mathcal{P}_{G\setminus e}$ with $w_{G\setminus
e}(\phi(f))=w_G(f)$.

(2) For any $f\in\mathcal{P}_G^0$, we have
$B_{G,\tau}(f)\setminus\{u\}=B_{G/e,\tau}(\phi(f))$.
\end{lem}

\begin{proof} (1) For any $I\subset V(G\setminus e)$ with $0\notin I$ and
$w\in I$, we have $outdeg_{I,G\setminus e}(w)=outdeg_{I,G}(w)$. This
implies $g=\phi(f)$ is a $(G\setminus e)$-parking function.
Conversely, for any $g\in\mathcal{P}_{G\setminus e}$, let
$f=\phi^{-1}(g)$ such that $f(w)=g(w)$ for any $w\in V(G\setminus
e)$ and $f(u)=0$. For any $I\subset V(G)$ with $0\notin I$, if $u\in
I$, then $f(u)<outdeg_{I,G}(u)$ since $f(u)=0$ and
$outdeg_{I,G}(u)\geq 1$; otherwise, we have $outdeg_{I,G\setminus
e}(w)=outdeg_{I,G}(w)$ for all $w\in I$, this implies
$f(w)<outdeg_{I,G}(w)$ for some $w\in I$ since $g$ is a $(G\setminus
e)$-parking function. Clearly, $w_{G\setminus e}(g)=w_G(f)$.

(2) Since $\tau(u)<\tau(w)$ for all $w\in NB_G(0)\setminus\{u\}$ and
$f(u)=0$, we have $Ord_{G,\tau}(f)_1=u$ for all
$f\in\mathcal{P}_G^0$. Let $g=\phi(f)$. Then $Ord_{G\setminus e,
\tau}(g)_i=Ord_{G, \tau}(f)_{i+1}$ for all $i\in[n-1]$. Let $v\in
B_{G,\tau}(f)\setminus\{u\}$. Clearly, $outdeg_{I_v,G\setminus
e}(v)=outdeg_{I_v,G}(v)$ and $f(v)=g(v)$. This implies
$g(v)=outdeg_{I_v,G\setminus e}(v)-1$. Hence, $v$
 is $g$-critical in $G\setminus e$ since $v$
 is $f$-critical in $G$.

Now, we assume that $W_{G\setminus e, v,\tau,g}\neq S_{G\setminus
e,v,\tau,g}$, i.e., there is a $g_1\in\mathcal{P}_{G\setminus e}$
such that $g_1\in W_{G\setminus e,v,\tau,g}$ and $g_1\notin
S_{G\setminus e,v,\tau,g}$. Suppose $Ord_{G,\tau}(f)_i=v$. So,
$Ord_{G\setminus e,\tau}(g)_i\neq v$. Furthermore, we have
$\tau(Ord_{G\setminus e,\tau}(g)_i)>\tau(v)$. Let
 $f_1=\phi^{-1}(g_1)$. It is easy to see that
$f_1\in W_{G,v,\tau,f}$ and $f_1\notin S_{G,v,\tau,f}$, a
contradiction.

Conversely, let $v\in B_{G\setminus e,\tau}(g)$. Clearly, $v\neq u$.
Assume $W_{G,v,\tau,f}\neq S_{G,v,\tau,f}$, where $f=\phi^{-1}(g)$,
i.e., there is a $f_1\in\mathcal{P}_G$ such that $f_1\in
W_{G,v,\tau,f}$ and $f_1\notin S_{G,v,\tau,f}$. Let $g_1=\phi(f_1)$.
Similarly, we can obtain $W_{G\setminus e,v,\tau,g}\neq
S_{G\setminus e,v,\tau,g}$, a contradiction.
\end{proof}

\begin{lem}\label{lemmabijectionu=1} (1) The mapping $\varphi$ is a bijection from $\mathcal{P}^1_G$
to $\mathcal{P}_{G-e}$ with $w_{G-e}(\varphi(f))=w_G(f)$.

(2) For any $f\in\mathcal{P}_G^1$, we have
$B_{G,\tau}(f)=B_{G\setminus e,\tau}(\varphi(f))$.

\end{lem}
\begin{proof} Since $f(u)\geq 1$, we have the edge $\{0,u\}$ isn't  a bridge. So,
  $G-e$ is still a connected graph. (1) For any $I\subset V(G-e)$ with $0\notin I$ and $w\in I$, we
have $outdeg_{I,G-e}(w)=outdeg_{I,G}(w)$ if $w\neq u$;
$outdeg_{I,G\setminus e}(w)=outdeg_{I,G}(w)-1$ if $w=u$. Note that
$g(u)=f(u)-1$. Hence, $g=\varphi(f)$ is a $(G-e)$-parking function.
Conversely, for any $g\in\mathcal{P}_{G-e}$, let $f=\varphi^{-1}(g)$
such that $f(w)=g(w)$ for any $w\in V(G-e)$ and $f(u)=g(u)+1$. For
any $I\subset V(G)$ with $0\notin I$, we have
$outdeg_{I,G}(w)=outdeg_{I,G-e}(w)$ if $w\neq u$;
$outdeg_{I,G}(w)=outdeg_{I,G-e}(w)+1$ if $w=u$. Note that
$f(u)=g(u)+1$. This implies $f(w)<outdeg_{I,G}(w)$ for some $w\in I$
since $g$ is a $(G-e)$-parking function. Clearly,
$w_{G-e}(g)=w_G(f)$.

(2) Note that $Ord_{G,\tau}(f)=Ord_{G-e,\tau}(\varphi(f))$ for all
$f\in\mathcal{P}_G$. For any $f\in\mathcal{P}_G$, it is easy to see
that the vertex $v$ is $\varphi(f)$-critical in $G-e$ if and only if
it is $f$-critical in $G$.

Now, given $f\in\mathcal{P}_G$, let $g=\varphi(f)$. For any $v\in
B_{G,\tau}(f)$, we assume that $W_{G-e,v,\tau,g}\neq
S_{G-e,v,\tau,g}$, i.e., there is a $g_1\in W_{G-e,v,\tau,g}$ and
$g_1\notin S_{G-e,v,\tau,g}$. Suppose $Ord_{G-e,\tau}(g)_i=v$. Then
$Ord_{G-e,\tau}(g_1)_i\neq v$. Furthermore, we have
$\tau(Ord_{G-e,\tau}(g_1)_i)>\tau(v)$. Let
 $f_1=\varphi^{-1}(g_1)$. It is easy to see that
$f_1\in W_{G,v,\tau,g}$  and $f_1\notin S_{G,v,\tau,g}$, a
contradiction.

Conversely, let $v\in B_{G-e,\tau}(g)$. Assume $W_{G,v,\tau,f}\neq
S_{G,v,\tau,f}$,
 i.e., there is a
$f_1\in\mathcal{P}_G$ such that $f_1\in W_{G,v,\tau,f}$ and
$f_1\notin S_{G,v,\tau,f}$. Let $g_1=\phi(f_1)$. Similarly, we can
obtain $W_{G-e,v,\tau,g}\neq S_{G-e,v,\tau,g}$, a contradiction.
\end{proof}

We are in a position to prove the main theorem.

\begin{thm}\label{theorem}Suppose that $G$ is a connected graph with vertex $\{0,1,\ldots,n\}$
 and $\tau$ is a vertex ranking in $S_n$. Let $T_G(x,y)$ be the Tutte polynomial of $G$. Then
$T_G(x,y)=\sum\limits_{f\in\mathcal{P}}x^{b(f)}y^{w(f)}$.
\end{thm}

\begin{proof} Let $P_{G}(x,y)=\sum\limits_{f\in\mathcal{P}}x^{b(f)}y^{w(f)}$. Let $e$ be an edge of $G$ connecting the vertices
 $u$ to $0$, where $u$ satisfies $\tau(u)\leq\tau(w)$ for all $w\in NB_G(0)$. We
consider the following three cases.\\
 {\it Case 1.} $e$ is a loop of $G$.

For any $f\in\mathcal{P}_G$, it is easy to see that $f$ is
$(G-e)$-parking function as well. Note that
\begin{eqnarray*}w_{G-e}(f)&=&|E(G-e)|-|V(G-e)|-\sum\limits_{i=0}^nf(i)\\
&=&|E(G)|-1-|V(G)|-\sum\limits_{i=0}^nf(i)\\
&=&w_{G}(f)-1.\end{eqnarray*} Hence,
\begin{eqnarray*}P_{G}(x,y)&=&\sum\limits_{f\in\mathcal{P}_G}x^{b(f)}y^{w(f)}\\
&=&\sum\limits_{g\in\mathcal{P}_{G-e}}x^{b_\tau(g)}y^{w(g)+1}\\
&=&yP_{G-e}(x,y)\end{eqnarray*}  {\it Case 2.} $e$ is a bridge of
$G$.

For any $f\in\mathcal{P}_G$, we have $f(u)=0$ since $e$ is a bridge.
So, $\mathcal{P}_G=\mathcal{P}_G^0$. Let $\phi$ be defined as that
in Lemma \ref{lemmabijectionu=0}. From Lemma \ref{nobridge}, we have
$u\in B_{G,\tau}(f)$ for all $f\in\mathcal{P}_G$. Lemma
\ref{lemmabijectionu=0} (2) tells us that $b_{G}(f)=b_{G\setminus
e}(\phi(f))+1$.

\begin{eqnarray*}P_{G}(x,y)&=&\sum\limits_{f\in\mathcal{P}_G}x^{b(f)}y^{w(f)}\\
&=&\sum\limits_{g\in\mathcal{P}_{G\setminus e}}x^{b(g)+1}y^{w(g)}\\
&=&xP_{G\setminus e}(x,y)\end{eqnarray*} {\it Case 3.} $e$ is
neither loop nor bridge of $G$.

First, we claim $u\notin B_G(f)$ for any $f\in\mathcal{P}_G^0$.
Since $\tau(u)\leq\tau(w)$ for all $w\in NB_G(0)$ and $f(u)=0$, we
have $Ord_{G,\tau}(f)_1=u$ for all $f\in\mathcal{P}_G^0$. If there
are at least two edges connecting the vertices $0$ and $u$ in $G$,
i.e., $\mu_G(0,u)\geq 2$,  then $u$ isn't $f$-critical for any
$f\in\mathcal{P}_G^0$. So, we suppose $\mu_G(0,u)=1$. This implies
that $u$ is $f$-critical for any $f\in\mathcal{P}_G^0$.  Since $e$
is neither loop nor bridge of $G$, there exists a
$f'\in\mathcal{P}_G$ such that $f'(u)\geq 1$. Suppose
$Ord_{G,\tau}(f')_1=w$. Then $w\neq u$, $f'(w)=0$ and
$\tau(w)>\tau(u)$. Hence,  $f'\in W_{G,u,\tau,f}$ and $f'\notin
S_{G,u,\tau, f}$. This tells us $u\notin B_{G}(f)$.

By Lemmas \ref{lemmabijectionu=0}  and \ref{lemmabijectionu=1}, we
have
\begin{eqnarray*}P_{G}(x,y)&=&\sum\limits_{f\in\mathcal{P}_G}x^{b(f)}y^{w(f)}\\
&=&\sum\limits_{f\in\mathcal{P}_{G}^0}x^{b(f)}y^{w(f)}+\sum\limits_{f\in\mathcal{P}_{G}^1}x^{b(f)}y^{w(f)}\\
&=&\sum\limits_{g\in\mathcal{P}_{G\setminus e}}x^{b(g)}y^{w(g)}+\sum\limits_{g\in\mathcal{P}_{G-e}}x^{b(g)}y^{w(g)}\\
&=&P_{G\setminus e}(x,y)+P_{G-e}(x,y)\end{eqnarray*}

Finally, we consider the initial conditions. Let $G$ be a graph with
vertex set $\{0\}$ and $E(G)=\emptyset$. There is a unique
$G$-parking function $f(0)=-1$. Clearly, $w(f)=0$ and
$B_G(f)=\emptyset$. So, $T_G(x,y)=1$. Next, let $G$ be a graph with
vertex set $\{0,1\}$ and $E(G)=\{\{0,1\}\}$. There is a unique
$G$-parking function $f(0)=-1$ and $f(1)=0$. It is easy to see
$B_G(f)=\{1\}$ and $w(f)=0$. Hence, $P_G(x,y)=x$. This complete the
proof.
\end{proof}
Let us define a multiset $BW_G=BW_{G,\tau}$ as
$BW_{G,\tau}=\{(b_{G,\tau}(f),w_{G}(f))\mid
f\in\mathcal{P}_\mathcal{G}\}.$ By Theorem \ref{theorem}, we
immediately obtain the following corollary.

\begin{cor} Let $G$ be a connected graph. Suppose $\tau_1$ and
$\tau_2$ are two vertex ranking. Then $BW_{G,\tau_1}=BW_{G,\tau_2}$.
\end{cor}

 Next, we consider the case in which $G$ is the complete graph
$K_{n+1}$ and $\tau$ is the identity permutation. Recall that the
$K_{n+1}$-parking functions are exactly the classical parking
functions, i.e., $f=(f(0),f(1),\ldots,f(n))\in\mathcal{P}_{K_{n+1}}$
if and only if $(f(1),\ldots,f(n))$ is a classical parking function.

\begin{defn}\label{classical} Given a classical parking function
$\alpha=(a_1,a_2,\ldots,a_n)$, we say that a term $a_i=j$ is
$\alpha$-critical maxima if $a_i$ satisfies that there are exactly
$n-1-j$ terms larger than $j$ and $k<i$ for all $a_k>j$.
\end{defn}

\begin{lem}  Let $\alpha=(a_1,a_2,\ldots,a_n)$ be a classical parking
function and $a_i=j$ an $\alpha$-critical maxima. Then there are
exactly $j$ terms less than $j$.
\end{lem}
\begin{proof} Let $\bar{\alpha}=(-1,a_1,a_2,\ldots,a_n)$. Then
$\bar{\alpha}\in\mathcal{P}_{K_{n+1}}$.  Let $T=\Phi(\bar{\alpha})$
be the spanning tree obtained by Algorithm A and $v=pre_T(i)$. Then
$Ord^{-1}(\bar{\alpha})_{v}=j$. Thus,
$Ord^{-1}(\bar{\alpha})_{i}=j+1$ since $a_i=j$ is an
$\alpha$-critical maxima. From Algorithm B,
$a_{Ord(\bar{\alpha})_k}<j$ for all $1\leq k\leq j$.
\end{proof}

\begin{lem} Let $\alpha=(a_1,a_2,\ldots,a_n)$ be a classical parking
function and $\bar{\alpha}=(-1,a_1,a_2,\ldots,a_n)$. Then $a_i=j$ is
$\alpha$-critical maxima if and only if the vertex $i$ is
$\bar{\alpha}$-bridge.
\end{lem}
\begin{proof} First, we suppose $a_i=j$ is
$\alpha$-critical maxima. By Algorithm A, it is easy to see the
vertex $i$ is $\bar{\alpha}$-critical since there are exactly $j$
terms less than $j$ and exactly $n-1-j$ terms larger than $j$. $k<i$
for all $a_k>a_i$ imply $W_{i,\bar{\alpha }}= S_{i,\bar{\alpha }}$.
Hence, the vertex $i$ is $\bar{\alpha}$-bridge.

Conversely, we suppose the vertex $i$ is $\bar{\alpha}$-bridge,
$a_i=j$ and $Ord^{-1}(\bar{\alpha})$.  Let $T=\Phi(\bar{\alpha})$ be
the spanning tree obtained by Algorithm A and $v=pre_T(i)$. There
are exactly $j$ terms less than $j$ and
$Ord^{-1}(\bar{\alpha})_v+1=Ord^{-1}(\bar{\alpha})_i$ since the
vertex $i$ is $\bar{\alpha}$-critical. Let $w$ be the first vertex
at the right of $i$
 such that $w>i$ in the order $Ord(\bar{\alpha})$. Let $T'$ be a new spanning tree from $T$ by deleting the
edge $\{v,i\}$ from $T$ and adding the edge $\{i,w\}$ into $T$ If
$a_w\leq j$; otherwise adding the edge $\{v,w\}$ into $T$. Let
$\bar{\beta}=\Theta(T')$ be the $K_{n+1}$-parking function by
Algorithm B. Then $\bar{\beta}\in W_{i,\bar{\alpha }}$ and
$\bar{\beta}\notin S_{i,\bar{\alpha }}$, a contradiction. So, we
prove that for any vertex $w\in [n]$ if $w$ is at the right of $i$
in the order $Org(\bar{\alpha})$, then $w<i$ and $a_w>j$. Hence,
$a_i$ is $\alpha$-critical maxima.
\end{proof}

Let $cm(\alpha)$ be the number of critical maxima in a classical
parking function $\alpha$.

\begin{cor}
\begin{eqnarray*}T_{K_{n+1}}(x,y)=\sum\limits_{\alpha\in\mathcal{P}_n}x^{cm(\alpha)}y^{{n\choose{2}}-\sum\limits_{i=1}^na_i}
\end{eqnarray*}
where $\mathcal{P}_n$ is the set of classical parking function of
length $n$.
\end{cor}
\begin{exa} We list all classical the parkin functions $\alpha$ of
length $3$ as well as the corresponding  sets of critical maxima and
$cm(\alpha)$ in the following table.
$$\begin{array}{|l|l|l|l|l|l|}
\hline
parking~function& critical~maxima & cm(\alpha)&parking~function& critical~maxima& cm(\alpha)\\
\hline (0,0,0)&\emptyset&0&(0,0,1)&\emptyset&0\\
\hline(0,0,2)&\{a_3\}&1&(0,1,0)&\emptyset&0\\
\hline(0,1,1)&\emptyset&0&(0,1,2)&\{a_3\}&1\\
\hline(0,2,0)&\{a_2\}&1&(0,2,1)&\{a_2,a_3\}&2\\
\hline(1,0,0)&\emptyset&0&(1,0,1)&\emptyset&0\\
\hline (1,0,2)&\{a_3\}&1&(1,1,0)&\{a_3\}&1\\
\hline (1,2,0)&\{a_2,a_3\}&2& (2,0,0)&\{a_1\}&1\\
\hline (2,0,1)&\{a_1,a_3\}&2& (2,1,0)&\{a_1,a_2,a_3\}&3\\
\hline
\end{array}
$$
\begin{center}
Table 4. Classical the parking functions $\alpha$ of length $n$ as
well as their sets of critical maxima
\end{center}
Hence, $T_{K_4}(x,y)=y^3+3y^2+2y+(4y+2)x+3x^2+x^3$.
\end{exa}

%%%%%%%%%%%%%%%%%%%%%%%%%%%%%%%%%%%%%%%%%%%%%%%%%%%%%%%%%%%%%%%%%%%%%%%%%%%%%%%


\begin{thebibliography}{99}

\bibitem{cori2002} R. Cori, D. Rossin, B. Salvy, Polynomial ideals
for sandpiles and their Grobner bases. {\it Theoretical Computer
Science} {\bf 276} (2002), no. 1-2, 1-15.



\bibitem{Denis2005} Chebikin Denis, Pylyavskyy Pavlo, A family of
bijections between $G$-parking functions and spanning trees. {\it J.
Combin. Theory Ser. A} {\bf 110} (2005), no. 1, 31--41.



\bibitem{kostic}Kosti$\acute{c}$ Dimitrije, Yan, Catherine H.
Multiparking functions, graph searching, and the Tutte polynomial.
{\it Adv. in Appl. Math.} {\bf 40} (2008), no. 1, 73--97.

\bibitem{EFL}Sen-Peng Eu,Tung-Shan Fu, Symmetric parking functions
and related multigraphs, {\it private communication}.



\bibitem{FR} D. Foata, J. Riordan, Mappings of acyclic and parking functions,
{\it Aequationes Math.} 10 (1974) 10-22.

\bibitem{F}J. Fran\c{c}n, Acyclic and parking functions, {\it J. Combin. Theory Ser.
A} 18 (1975) 27-35.

\bibitem{GK}J.D. Gilbey, L.H. Kalikow, Parking functions, valet functions and
priority queues, {\it Discrete Math.} 197/198 (1999) 351-373.




\bibitem{konhein1966} Konheim, A. G. and Weiss, B. An Ocuupancy
Discipline and Applications. {\it Siam Journal of Applied
Mathematics} {\bf 14} (1966) 1266-1274





\bibitem{KY}J.P.S. Kung, C.H. Yan, Gon\v{c}arove polynomials and parking functions,
{\it J. Combin. Theory Ser. A} 102 (2003) 16-37.

\bibitem{postnikov2004} Postnikov, A. and Shapiro, B. Trees, Parking
Functions, Syzygies, and Deformatioins of Monomial Ideals. {\it
Transactions of the American Mathematical Society} {\bf 356} (2004).

\bibitem{PS}J. Pitman, R. Stanley, A polytope related to empirical
distributions, plane trees, parking functions, and the
associahedron, {\it Discrete Comput. Geom.} 27 (4) (2002) 603-634.

\bibitem{R} J.Riordan, Ballots and trees, {\it J.Combin. Theory} 6 (1969) 408-411.

\bibitem{SMP}M.P. Sch¨¹tzenberger, On an enumeration problem, {\it J. Combin. Theory} 4 (1968) 219-221.

\bibitem{SRP}R.P. Stanley, Hyperplane arrangements, interval orders and trees,
{\it Proc. Natl. Acad. Sci.} 93 (1996) 2620-2625.

\bibitem{SRP2}R.P. Stanley, Parking functions and non-crossing partitions, in: The
Wilf Festschrift, {\it Electron. J. Combin.} 4 (1997) R20.



\bibitem{Y1} C.H. Yan, Generalized tree inversions and k-parking functions, {\it J.
Combin. Theory Ser. A} 79 (1997) 268-280.

\bibitem{Y2}C.H. Yan, On the enumeration of generalized parking functions,
{\it Congr. Numer.} 147 (2000) 201-209.

\bibitem{Y3}C.H. Yan, Generalized parking functions, tree inversions and
multicolored graphs, {\it Adv. in Appl. Math.} 27 (2001) 641-670.

\end{thebibliography}
\end{document}